\documentclass[11pt]{article}
\usepackage{amsmath,amsfonts,amsthm,amssymb,setspace}
\newcommand\Tr{{\rm Tr\,}}
\newtheorem{thm}{Theorem}
\newtheorem{rem}{Remark}
\begin{document}
\newcommand{\email}[1]{{\sl E-mail:\/} {\texttt{#1}}}
\newcommand{\version}{July 24, 2007}

%%%%%%%%%%%%%%%%%%%%%%%%%%%%%%%%%%%%%%%%%%%%%%%%%%%%%%%%%%%%%%%%%%%%%%%%%%%%%%%%
\title{Lieb-Thirring inequalities with improved constants}
\author{Jean Dolbeault, Ari Laptev, and Michael Loss}
\date{}
\maketitle
\thispagestyle{empty}
\begin{abstract} Following Eden and Foias we obtain a matrix version of a generalised Sobolev inequality in one-dimension. This allows us to improve on the known estimates of best constants in Lieb-Thirring inequalities for the sum of the negative eigenvalues for multi-dimensional Schr\"odinger operators.
\end{abstract}

\noindent\begin{minipage}{125mm}\linespread{0.9}\selectfont{\small {\sl Key-words:\/} Sobolev inequalities; Schr\"odinger operator; Lieb-Thirring inequalities. {\sl MSC (2000):\/} Primary: 35P15; Secondary: 81Q10}\end{minipage}

%%%%%%%%%%%%%%%%%%%%%%%%%%%%%%%%%%%%%%%%%%%%%%%%%%%%%%%%%%%%%%%%%%%%%%%%%%%%%%%%
\section{Introduction}

Let $H$ be a Schr\"odinger operator in $L^2(\mathbb{R}^d)$
\begin{equation}\label{H}
H = -\Delta - V
\end{equation}
For a real-valued potential $V$ we consider Lieb-Thirring inequalities
for the negative eigenvalues $\{\lambda_n\}$ of the operator $H$
\begin{equation}\label{LTh}
\sum |\lambda_n|^\gamma \le L_{d,\gamma} \int_{\mathbb{R}^d}
V_+^{d/2+\gamma}(x) \, dx\,,
\end{equation}
where $V_+ = (|V| + V)/2$ is the positive part of $V$.

Eden and Foias have obtained in \cite{EF} a version of a one-dimensional generalised Sobolev inequality which gives best known estimates for the constants in the inequality \eqref{LTh} for $1\le\gamma<3/2$. The aim of this short article is to extend the method from \cite{EF} to a class of matrix-valued potentials. By using ideas from \cite{LW} this automatically improves on the known estimates of best constants in \eqref{LTh} for multidimensional Schr\"odinger operators.

Lieb-Thirring inequalities for matrix-valued potentials for the value $\gamma=3/2$ were obtained in \cite{LW} and also in \cite{BL}. Here we state a result corresponding to $\gamma=1$.

%-------------------------------------------------------------------------------
\begin{thm}\label{1D-matrix-LTh}
Let $V\ge 0$ be a Hermitian $m\times m$ matrix-function defined on~$\mathbb R$ and let $\lambda_n$ be all negative eigenvalues of the operator \eqref{H}. Then
\begin{equation}\label{Th1}
\sum |\lambda_n| \le \frac{2}{3\sqrt 3} \int_{\mathbb R}\Tr\left[V^{3/2}(x)\right] \, dx\,.
\end{equation}
\end{thm}
%-------------------------------------------------------------------------------

\begin{rem}
The constant $\frac 2{3\sqrt3}$ should be compared with the Lieb-Thirring constant found in \cite{LT} for a class of single eigenvalue potentials and with the constant obtained in \cite{HLW} which is twice as large as the semi-classical one
\begin{equation*}
\frac 4{3\sqrt 3\,\pi}<\frac{2}{3\sqrt 3}<2\times \frac{2}{3\pi} = 2\times \frac 1{2\pi} \int_{\mathbb R} (1-\xi^2)_+\,d\xi\,.
\end{equation*}
This is about $0,2450\dots<0,3849\dots <0,4244\dots$.
\end{rem}

\begin{rem}
Note that the values of the best constants for the range $1/2<\gamma<3/2$ remain unknown.
\end{rem}

Let $\mathcal A(x) = (a_1(x),\dots, a_d(x))$ be a magnetic vector potential with real valued entries $a_k\in L^2_{\rm loc}(\mathbb R^d)$ and let
\begin{equation*}\label{H(A)}
H(\mathcal A) = (i\,\nabla + \mathcal A)^2 -V\,,
\end{equation*}
where $V\ge 0 $ is a real-valued function.

Denote the ratio of $2/3\sqrt3$ and the semi-classical constant by
\begin{equation*}
R := \frac{2}{3\sqrt 3} \times \left(\frac 2{3\pi}\right)^{-1} = 1.8138 \dots.
\end{equation*}

By using the Aizenmann-Lieb argument \cite{AL}, a ``lifting" with respect to dimension \cite{LW}, \cite{HLW}, and Theorem \ref{1D-matrix-LTh} we obtain the following result:

%-------------------------------------------------------------------------------
\begin{thm}\label{d-LTh}
For any $\gamma \ge 1$ and any dimension $d \ge 1$, the negative eigenvalues of the operator $H(\mathcal A)$ satisfy inequalities
\begin{equation*}
\sum |\lambda_n|^\gamma \le L_{d,\gamma} \int_{\mathbb R^d} V^{d/2 + \gamma}(x) \, dx\,,
\end{equation*}
where
\begin{equation*}
L_{d,\gamma} \le R\times L_{d,\gamma}^{cl} = R\times \frac 1{(2\pi)^d} \int_{\mathbb R^d} (1-|\xi|)_+^\gamma\, d\xi\,.
\end{equation*}
\end{thm}
%-------------------------------------------------------------------------------

\begin{rem}
Theorem \ref{d-LTh} allows us to improve on the estimates of best constants in Lieb-Thirring inequalities for Schr\"odinger operators with complex-valued potentials recently obtained in \cite{FLLS}.
\end{rem}

%%%%%%%%%%%%%%%%%%%%%%%%%%%%%%%%%%%%%%%%%%%%%%%%%%%%%%%%%%%%%%%%%%%%%%%%%%%%%%%%
\section{One-dimensional generalised Sobolev inequality for matrices}

Let $\{\phi_n\}_{n=0}^N$ be an ortho-normal system of vector-functions in $L^2(\mathbb R, \mathbb C^M)$, $M\in\mathbb N$,
\begin{equation*}
(\phi_n,\phi_m) = (\phi_n,\phi_m)_{L^2(\mathbb R, \mathbb C^M)} = \sum_{j=1}^M \int_{\mathbb R} \phi_n(x,j)\,\overline{\phi_m(x,j)}\, dx = \delta_{nm}\,,
\end{equation*}
where $\delta_{nm}$ is the Kronecker symbol. Let us introduce an $M\times M$ matrix $U$ with entries
\begin{equation*}
u_{j,k}(x,y) = \sum_{n=0}^N \phi_n(x,j)\,\overline{\phi_n(y,k)}\,.
\end{equation*}
Clearly
\begin{equation}\label{P*=P}
U^*(x,y) = U(y,x)\,.
\end{equation}
The fact that the functions $\phi_n$ are orthonormal can be written in a compact form
\begin{equation}\label{PP=P}
\int_{\mathbb R} U(x,y)\,U(y,z)\, dy = U(x,y)\,.
\end{equation}
The latter two properties \eqref{P*=P} and \eqref{PP=P} prove that $U(x,y)$ could be considered as the integral kernel of an orthogonal projection $P$ in $L^2(\mathbb R, \mathbb C^M)$ whose image is the subspace of vector-functions spanned by $\{\phi_n\}_{n=1}^N$.

%-------------------------------------------------------------------------------
\begin{thm}\label{1D-Sobolev}
Let us assume that the vector-function $\phi_n$, $n=1,2,\dots N,$ are from the Sobolev class $H^1(\mathbb R, \mathbb C^M)$. Then
\begin{equation*}\label{1d-sobolev}
\int_{\mathbb R} \Tr\Big[U(x,x)^3\Big]\,dx \le \sum_{n=1}^N \sum_{j=1}^M\, \int_{\mathbb R} |\phi_n'(x,j)|^2\,dx\,.
\end{equation*}
\end{thm}
%-------------------------------------------------------------------------------

\begin{proof}
\begin{multline}\label{UUU}
\frac{d}{dy}\, \Tr\Big[U(x,y)\,U(y,x)\, U(x,x)\Big]\\ =
\Tr\! \left[\Big(\frac{d}{dy}\, U(x,y)\Big)\,U(y,x)\, U(x,x)\right] + \Tr\! \left[U(x,y)\,\Big( \frac{d}{dy}\, U(y,x)\Big)\,U(x,x)\right]
\end{multline}
By integrating \eqref{UUU} and taking absolute values one obtains
\begin{multline*}
\frac 12\,\Tr\Big[ U(x,z)\,U(z,x)\,U(x,x)\Big]\\ \le
\frac12\, \int_{-\infty}^z \Big |\,\Tr\Big[ \Big(\frac{d}{dy}\, U(x,y)\Big)\,U(y,x)\, U(x,x)\Big]\\
+ \Tr\Big[U(x,y)\,\Big( \frac{d}{dy}\, U(y,x)\Big)\,U(x,x)\Big] \Big |\, dy
\end{multline*}
and
\begin{multline*}
\frac 12\,\Tr\Big[ U(x,z)\,U(z,x)\,U(x,x)\Big]\\ \le
\frac12\,\int_z^\infty \Big|\, \Tr\Big[ \Big(\frac{d}{dy}\, U(x,y)\Big)\,U(y,x)\, U(x,x)\Big]\\
+ \Tr\Big[ U(x,y)\,\Big( \frac{d}{dy}\, U(y,x)\Big)\,U(x,x)\Big]\Big|\, dy\,.
\end{multline*}
Taking absolute values and adding the two inequalities yields for any
$z \in \mathbb{R}$
\begin{multline} \label{int++UUU}
\Big|\,\Tr\Big[ U(x,z)\,U(z,x)\,U(x,x)\Big]\Big|\\
\le \frac12 \, \int_{\mathbb R} \left|\,\Tr\left[\left(\frac{d}{dy}\, U(x,y)\right) U(y,x)\, U(x,x)\right]\right|\, dy\\
+ \frac12 \, \int_{\mathbb R} \left|\,\Tr\left[ U(x,y) \left(\frac{d}{dy}\,U(y,x) \right) U(x,x)\right]\right|\, dy\,.
\end{multline}
Note that we have reproved Agmon's inequality
$$
|f(x)|^2 \le \int_{\mathbb R} |f(y)\,f'(y)|\, dy
$$
for traces of matrices. By using properties of traces, the Cauchy-Schwarz inequality for matrix-functions and also properties \eqref{P*=P} and \eqref{PP=P}, we find that for all $z \in \mathbb{R}$
\begin{multline*}
\left(\int_{\mathbb R} \left|\,\Tr\left[ \left(\frac{d}{dy}\, U(x,y)\right)\,U(y,x)\,U(x,x)\right]\right|\, dy\right)^2\\
\le \int_{\mathbb R} \Tr\left[\frac{d}{dy}\, U(x,y)^*\, \frac{d}{dy}\, U(x,y)\right] dy\; \int_{\mathbb R} \Tr \left[U(x,y)^*\,U^2(x,x)\,U(x,y)\right] dy\\
= \int_{\mathbb R} \Tr \left[\frac{d}{dy}\, U(y,x)\, \frac{d}{dy}\, U(x,y)\right] dy\; \int_{\mathbb R} \Tr \left[U^2(x,x)\,U(x,y)\,U(y,x)\right] dy\\
=\int_{\mathbb R} \Tr \left[\frac{d}{dy}\, U(x,y)\,\frac{d}{dy}\, U(y,x)\right] dy\;\;\Tr \left[U(x,x)^3\right] ,
\end{multline*}
and similarly
\begin{multline*}
\left(\int_{\mathbb R} \left|\,\Tr\left[ U(x,y)\,\frac{d}{dy}\,U(y,x)\,U(x,x)\right]\right|\, dy\right)^2\\
\le \int_{\mathbb R} \Tr \left[\frac{d}{dy}\,U(x,y)\, \frac{d}{dy}\,U(y,x)\right] dy\;\;\Tr \left[U(x,x)^3\right]\,.
\end{multline*}
Thus, using this, and setting $x=z$ in \eqref{int++UUU}, we arrive at
\begin{equation*}
\Big|\,\Tr\Big[ U(x,x)^3\Big] \Big| \le \int_{\mathbb R} \Tr \left[\frac{d}{dy}\,U(x,y)\, \frac{d}{dy}\,U(y,x)\right]\, dy \,.
\end{equation*}
Integrating with respect to $x$ we finally obtain
\begin{multline*}
\int_{\mathbb R} \Big|\,\Tr\Big[ U(x,x)^3\Big] \Big|\,dx\\
\le \sum_{n,k=1}^N\sum_{i,j=1}^M \int_{\mathbb R} \int_{\mathbb R} \phi_n(x,i)\, \overline{\phi_n'(y,j)}\; \phi_k'(y,j)\, \overline{\phi_k(x,i)} \, dx\,dy\\
= \sum_{n=1}^N\sum_{j=1}^M \int_{\mathbb R} |\phi_n'(x,j)|^2\, dx\,,
\end{multline*}
which completes the proof.
\end{proof}

%%%%%%%%%%%%%%%%%%%%%%%%%%%%%%%%%%%%%%%%%%%%%%%%%%%%%%%%%%%%%%%%%%%%%%%%%%%%%%%%
\section{Lieb-Thirring inequalities for Schr\"odinger operators with matrix-valued potentials}

Let us assume that $V\in C_0^\infty(\mathbb R)$, $V\ge0$, be a $M\times M$ Hermitian matrix-valued potential with entries $\{v_{ij}\}_{i,j=1}^M$. Then the negative spectrum of the Schr\"odinger operator $H = -\frac{d^2}{dx^2} - V$ in $L^2(\mathbb R)$ is finite.

Denote by $\{\phi_n\}$ the ortho-normal system of eigen-vector functions corresponding to the eigenvalues $\{\lambda_n\}_{n=1}^N$
\begin{equation*}
-\frac{d^2}{dx^2}\,\phi_n - V\phi_n = \lambda_n\,\phi_n\,.
\end{equation*}
Clearly,
\begin{equation*}
\sum_n \lambda_n = \sum_{n,j} \int_{\mathbb R} |\phi_n'(x,j)|^2 \, dx - \Tr\left[ \int_{\mathbb R}V(x)\,U(x,x)\, dx\right]
\end{equation*}
and by H\"older's inequality for traces,
\begin{equation*}
\int_{\mathbb R} \Tr\left[ V(x)\,U(x,x)\right]\, dx \le \left(\int_{\mathbb R}\Tr\big[ V^{3/2}(x)\big]\, dx \right)^\frac 23 \!\left(\int_{\mathbb R} \Tr\left[ U(x,x)^3\right]\, dx\right)^\frac 13\!,
\end{equation*}
so that using Theorem \ref{1D-Sobolev}
\begin{equation*}
\sum_n \lambda_n \ge X - \left(\int_{\mathbb R}\Tr\left[ V^{3/2}(x)\right]\, dx \right)^\frac 23 \,X^\frac 13
\end{equation*}
with $ X := \int_{\mathbb R} \Tr\left[ U(x,x)^3\right] dx$. Minimising the right hand side with respect to $X$ we finally complete the proof of Theorem \ref{1D-matrix-LTh}
\begin{equation*}
\sum_n \lambda_n \ge -\frac{2}{3\sqrt3} \int_{\mathbb R}\Tr\left[V^{3/2}(x)\right]\, dx\,.
\end{equation*}

%%%%%%%%%%%%%%%%%%%%%%%%%%%%%%%%%%%%%%%%%%%%%%%%%%%%%%%%%%%%%%%%%%%%%%%%%%%%%%%%
\noindent\begin{minipage}{125mm}\linespread{0.9}\selectfont
{\it Acknowledgements.} {\small The authors are grateful to the organisers of the meeting \lq\lq Functional Inequalities: Probability and PDE's", Universit\'e Paris-X, June 4-6, 2007, where this paper came to fruition. We would like to thank Robert Seiringer for pointing out an omission in the formulation of Theorem 2. A.L. thanks the Department of Mathematics of the University Paris Dauphine for its hospitality and also the ESF Programme SPECT. M.~L. would like to acknowledge partial support through NSF grant DMS-0600037.}\end{minipage}

%%%%%%%%%%%%%%%%%%%%%%%%%%%%%%%%%%%%%%%%%%%%%%%%%%%%%%%%%%%%%%%%%%%%%%%%%%%%%%%%

\small
\noindent
{\sc J. Dolbeault:} Ceremade UMR CNRS no. 7534, 
Universit\'eŽ Paris Dauphine, F-75775 Paris Cedex 16, France.
\email{dolbeaul@ceremade.dauphine.fr}
\par\smallskip\noindent
{\sc A. Laptev:} Department of Mathematics,
Imperial College London, London SW7 2AZ, UK,
Royal Institute of Technology, 100 44 Stockholm, Sweden.
\email{a.laptev@imperial.ac.uk}
\par\smallskip\noindent
{\sc M. Loss:} School of Mathematics,
Georgia Institute of Technology, Atlanta GA 30332-0160, USA.
\email{loss@math.gatech.edu}

\begin{flushright}{\sl\version}\end{flushright}
\end{document}